\documentclass[11pt]{article}

\PassOptionsToPackage{superscript,biblabel}{cite}
\PassOptionsToPackage{hidelinks}{hyperref}

\usepackage[scaled=.98,p]{XCharter}
\usepackage[scaled=1.04,varqu,varl]{inconsolata}
\usepackage[type1]{cabin}
\usepackage[uprightscript,xcharter,vvarbb,scaled=1.05]{newtxmath}
\usepackage[T1]{fontenc}
\usepackage[final,protrusion=true,expansion=true]{microtype}

\usepackage{titling}
\usepackage{authblk}
\pretitle{\begin{center}\bfseries\sffamily\LARGE}
\posttitle{\end{center}}
\usepackage{abstract}

\usepackage{sectsty}
\allsectionsfont{\sffamily}

\usepackage[left=1in,right=1in,top=1in,bottom=1in,includefoot,heightrounded]{geometry}
\usepackage{ifpdf}
\newcommand{\CAVArXivMode}{}
\usepackage{hyperref}
\usepackage{bm}
\usepackage{booktabs}
\usepackage[margin=15pt,font=small,labelfont={bf,sf},justification=justified]{caption}
\usepackage{cite}
\usepackage{graphicx}
\usepackage{mathtools}
\ifdefined\CAVArXivMode
\else
  \usepackage[bbgreekl]{mathbbol}
\fi
\ifdefined\CAVArXivMode
\else
  \usepackage{amsfonts}
\fi
\usepackage{mathrsfs}
\usepackage{algorithm}
\usepackage{algorithmic}
\usepackage{tikz}
\usetikzlibrary{arrows.meta}
\usepackage{subcaption}
\usepackage[export]{adjustbox}
\usepackage{todonotes}

\ifpdf
  \DeclareGraphicsExtensions{.eps,.pdf,.png,.jpg}
\else
  \DeclareGraphicsExtensions{.eps}
\fi

\newcommand{\Ab}{\bm{\mathsf{A}}}
\newcommand{\Azero}{\Ab_0}
\newcommand{\Bb}{\bm{\mathsf{B}}}
\newcommand{\Cb}{\bm{\mathsf{C}}}
\newcommand{\Db}{\bm{\mathsf{D}}}
\newcommand{\Eb}{\bm{\mathsf{E}}}
\newcommand{\Ezero}{\Eb_0}
\newcommand{\Eeul}{\Eb_{\text{E}}}
\newcommand{\Eeulzero}{\Eb_{\text{E},0}}
\newcommand{\Fb}{\bm{F}}
\newcommand{\Ub}{\bm{U}}
\newcommand{\Gb}{\bm{\mathsf{G}}}
\newcommand{\Ib}{\bm{\mathsf{I}}}
\newcommand{\Jb}{\bm{\mathsf{J}}}
\newcommand{\Kb}{\bm{\mathsf{K}}}
\newcommand{\KE}{\Kb_{\text{E}}}
\newcommand{\Lb}{\bm{\mathsf{L}}}
\newcommand{\Mb}{\bm{\mathsf{M}}}
\newcommand{\Pb}{\bm{\mathsf{P}}}
\newcommand{\Ph}{\mathscr{P}_{h}^{0}}
\newcommand{\Rb}{\bm{\mathsf{R}}}
\newcommand{\Sb}{\bm{\mathsf{S}}}
\newcommand{\Tb}{\bm{\mathsf{T}}}
\newcommand{\Thetab}{\bm{\mathsf{\Theta}}}
\newcommand{\Uh}{\mathscr{U}_{h}}
\newcommand{\Dlink}{\bm{D}}
\newcommand{\bb}{\bm{b}}
\newcommand{\deltab}{\bm{\delta}}
\newcommand{\fb}{\bm{f}}
\newcommand{\rb}{\bm{r}}
\newcommand{\ub}{\bm{u}}
\newcommand{\vb}{\bm{v}}
\newcommand{\wb}{\bm{w}}
\newcommand{\xb}{\bm{x}}
\newcommand{\Xb}{\bm{X}}
\newcommand{\zerob}{\bm{0}}
\newcommand{\zeromat}{\bm{\mathsf{0}}}
\newcommand{\ds}{\Delta s}
\newcommand{\dt}{\Delta t}
\newcommand*{\adj}{^{*}}
\newcommand*{\tran}{^{\mkern-1.5mu\mathsf{T}}}
 \providecommand{\bbsigma}{\mathbb{\sigma}}
\DeclareSymbolFont{CAVsansGreek}{OML}{cmssm}{b}{it}
\DeclareMathSymbol{\CAVTheta}{\mathord}{CAVsansGreek}{"02}
\renewcommand{\Thetab}{\CAVTheta}
\usepackage{xurl}
\usepackage[capitalise,nameinlink,noabbrev]{cleveref}
\newcommand{\TheTitle}{Coupling-Aware Vanka Smoothing for Multigrid Preconditioning of the Implicit Immersed Boundary Equations
}

\newcommand{\TheFunding}{The work of the first author was supported by the Department of Defense (DoD) through the National Defense Science and Engineering Graduate (NDSEG) Fellowship Program. The work of the second author is supported by the National Science Foundation (grant numbers OAC 1450327, OAC 1652541, OAC 1931516, and DMS 2608488) and the National Institutes of Health (NIH) (grant numbers HL157631 and HL182166). This manuscript is the result of funding in whole or in part by the NIH. It is subject to the NIH Public Access Policy. Through acceptance of this federal funding, NIH has been given a right to make this manuscript publicly available in PubMed Central upon the Official Date of Publication, as defined by NIH.
}
 
\newenvironment{keywords}{\par\medskip\noindent\textbf{Keywords.}\ }{\par}
\newenvironment{MSCcodes}{\par\medskip\noindent\textbf{2020 Mathematics Subject Classification.}\ }{\par\medskip}
\title{\TheTitle}
\author[1]{Cole Gruninger}
\author[2]{Boyce E. Griffith}
\affil[1]{Department of Mathematics, Fordham University}
\affil[2]{Departments of Mathematics and Biomedical Engineering, University of North Carolina at Chapel Hill}
\date{}

\ifpdf
\hypersetup{
  pdftitle={\TheTitle},
  pdfauthor={C. Gruninger, B. E. Griffith}
}
\fi

\begin{document}

\maketitle
\begin{abstract}
The immersed boundary (IB) method models fluid--structure interaction using the natural Lagrangian and Eulerian formulations of structural mechanics and fluid dynamics, respectively, but explicit time discretization of the IB force imposes a stiffness-dependent upper bound on the time-step size. Treating these forces implicitly removes this restriction, but the resulting coupled linear systems become increasingly difficult to solve as the structural stiffness increases. Algebraically eliminating the Lagrangian degrees of freedom yields a reduced Eulerian velocity--pressure IB system. Here we introduce a coupling-aware Vanka (CAV) smoothing strategy to enable effective multigrid preconditioning of this system. CAV patches are built as unions of standard pressure-centered Vanka patches, with the graph of the Eulerian elasticity matrix determining which patches are combined. Under grid refinement, CAV patch sizes remain bounded, and the computational cost of each multigrid cycle thereby grows linearly with the number of Eulerian degrees of freedom. Tests using target-point, membrane, and beam force laws show that CAV-preconditioned FGMRES reduces the relative residual by ten orders of magnitude in \(9\)--\(15\) iterations, with little growth under grid refinement. In a nonlinear benchmark modeling flow past a flexible fiber, the average number of FGMRES iterations per Newton solve increases only from \(8.6\) to \(9.5\) as the Eulerian grid is refined from \(32\times32\) to \(256\times256\) cells. To our knowledge, CAV provides the first robust multigrid strategy for time-dependent implicit IB formulations.
\end{abstract}

\begin{keywords}
multigrid, patch relaxation, Vanka smoothing, immersed boundary method, fluid--structure interaction, saddle-point systems, preconditioning
\end{keywords}

\begin{MSCcodes}
65M55, 65F08, 65F10, 65M85, 74F10
\end{MSCcodes}

\section{Introduction}
\label{sec:intro}

Peskin's immersed boundary (IB) method~\cite{peskin1972,peskin2002} provides a framework for modeling fluid--structure interaction involving elastic structures and viscous, incompressible fluids. Explicit time discretizations are straightforward to implement, but the stiff nature of immersed boundary forces can impose a severe stability restriction on the time-step size~\cite{stockie1999,hua2022}. Fully implicit and semi-implicit schemes alleviate this limitation, allowing the time step to be chosen based on other considerations rather than the stiffness of the elastic forces~\cite{newren2007a,newren2008,hou2008,mori2008a,ceniceros2011,ceniceros2009}. For instance, with the semi-implicit fluid discretization used here, the convective CFL condition determines the time-step scaling. However, fully implicit discretizations are often challenging to solve because the arguments of the regularized delta function introduce nonlinear dependence. Semi-implicit schemes therefore offer an attractive compromise: they treat the arguments of the regularized delta function explicitly and discretize the Lagrangian force density and interpolated velocity implicitly. Newren et al.~\cite{newren2007,newren2007a} proved unconditional linear stability of these schemes, provided that the discrete fluid--structure coupling preserves power between the Eulerian and Lagrangian frames and that the force law defines a nonnegative elastic energy functional. Algebraically, these two requirements correspond to the IB force-spreading and velocity-interpolation operators being adjoints of one another, and to the Lagrangian force operator being self-adjoint and negative semidefinite in the quadrature-weighted Lagrangian inner product.

Although this semi-implicit approach removes the time-step stability constraint imposed by the elastic forces, its formulation leads to linear systems that become increasingly ill-conditioned as the stiffness of the immersed structure increases. Effective multigrid preconditioners for these systems are lacking, particularly in the time-dependent setting. Consequently, the cost of the linear solves can outweigh the benefit of using larger time steps.

Here, we develop an effective coupling-aware Vanka (CAV) multigrid preconditioner for an Eulerian Schur complement formulation of the semi-implicit IB equations~\cite{newren2007a}, in which the Lagrangian degrees of freedom have been algebraically eliminated. The Eulerian degrees of freedom, namely the values of the velocity and pressure, lie on a Cartesian marker-and-cell (MAC) grid, so the system retains the grid structure needed for geometric multigrid methods. At each implicit time step, the velocity--pressure system is governed by the reduced Eulerian operator
\begin{equation}
\label{eq:intro_saddle}
\KE
=
\begin{bmatrix}
\Ab & \Gb\\
-\Db & \zeromat
\end{bmatrix},
\qquad
\Ab=\Azero-\dt\,\Eeul,
\qquad
\Azero=\frac{\rho}{\dt}\Ib-\mu\Lb,
\qquad
\Eeul=\Sb\Eb\Jb.
\end{equation}
Here, \(\Db\), \(\Gb\), and \(\Lb\) are the discrete divergence, gradient, and vector Laplacian matrices, respectively; \(\Sb\) is the IB force-spreading matrix that maps Lagrangian force densities to Eulerian force densities; \(\Jb=\Sb\adj\) is the IB velocity-interpolation matrix that maps Eulerian velocities to Lagrangian marker velocities; and \(\Eb\) is the Lagrangian elasticity matrix, with \(\Fb=\Eb\Xb\) for a linear force law and \(\delta\Fb=\Eb\,\delta\Xb\) in the linearization of a nonlinear force law. The matrix \(\Azero\) is the time-discrete fluid momentum matrix in the absence of elasticity, and \(\Eeul\) represents the elastic response on the Eulerian grid. Their combination \(\Ab=\Azero-\dt\,\Eeul\) forms the velocity block of the reduced Eulerian saddle-point system. Because \(\Sb\) and \(\Jb\) are adjoints, if \(\Eb\) is self-adjoint and negative semidefinite in the quadrature-weighted Lagrangian inner product, then \(\Eeul\) is self-adjoint and negative semidefinite in the Eulerian inner product, and \(-\dt\,\Eeul\) is positive semidefinite; hence, \(\Ab\) is obtained by shifting \(\Azero\) by the positive-semidefinite matrix \(-\dt\,\Eeul\). Thus, \(\KE\) has a familiar Stokes-like saddle-point structure, but its momentum block also contains the sparse Eulerian elasticity term \(-\dt\,\Eeul\). At high structural stiffness, this term dominates \(\Ab\), making its treatment central to effective preconditioning.

The block structure of \(\KE\) naturally suggests projection and block-factorization preconditioners. Both approaches require an approximation to the pressure Schur complement \(\Db\Ab^{-1}\Gb\). For the MAC-grid fluid operator \(\Azero\), commutator-based approximations choose a pressure operator \(\Ab_{0,\mathrm p}\) satisfying \(\Azero\Gb\approx\Gb\Ab_{0,\mathrm p}\), with a commutator error independent of the grid resolution. Including elasticity leaves a commutator residual proportional to \(\Eeul\Gb\), whose magnitude increases linearly with the structural stiffness. IB-specific projection preconditioners can remain effective at low Reynolds number and high stiffness, but their reported efficiency deteriorates if both Reynolds number and stiffness are large~\cite{zhang2014}. Alternatively, Luhana and Greif~\cite{luhana2026} retain the Lagrangian variables and develop a triangular block preconditioner for the corresponding monolithic double saddle-point system. Their factorization isolates a Lagrangian Schur complement that represents the fluid-mediated response of the immersed structure. They approximate this response by replacing the full-domain Stokes solve with a solve on a rectangular subdomain of the original Cartesian grid surrounding the structure. The resulting method combines nested Krylov iterations with standard multigrid-preconditioned fluid solvers. Using the truncated domain approximation substantially reduces the cost of the nested fluid solves, and their refinement studies exhibit approximately linear growth in solve time with the number of degrees of freedom.

Multigrid methods have long provided highly effective solvers for Poisson-type problems, whereas correspondingly robust multigrid methods for saddle-point systems were developed more recently. Vertex-star relaxation, pioneered by Schöberl and Zulehner~\cite{schoberl2003saddle}, associates a patch with each mesh vertex and collects the degrees of freedom supported on the cells incident to that vertex, together with the corresponding constraint unknowns. Farrell et~al.~\cite{farrell2021pcpatch} systematized this topological construction in \texttt{PCPATCH}, which is implemented in PETSc~\cite{petsc-user-ref}.

Pressure-centered Vanka relaxation has an analogous interpretation on a MAC grid: a standard pressure-centered Vanka patch contains one pressure unknown together with the velocity unknowns in its discrete divergence equation~\cite{vanka1986,manservisi2006vanka}. Guy et~al.~\cite{guy2015} generalized this standard pressure-centered patch to larger geometric boxes for the implicit IB equations. This box-relaxation scheme captures more of the coupling introduced by \(\Eeul\), but because the boxes are fixed \emph{a priori}, their boundaries can still cut through \(\Eeul\) stencils. Consequently, as illustrated here, box relaxation does not yield a grid-scalable preconditioner for time-dependent IB models with membrane or beam force laws.

CAV patches are built as unions of standard pressure-centered Vanka patches, with the graph of \(\Eeul\) determining which patches are combined. Starting from each standard pressure-centered Vanka patch, CAV first augments its velocity set with every velocity coupled through a nonzero entry of \(\Eeul\) to a velocity already in the patch. If the augmentation introduces no additional velocity degrees of freedom, CAV retains only the standard pressure-centered Vanka patch. Otherwise, it identifies the pressures whose discrete divergence equations involve the expanded velocity set and takes the union of the standard pressure-centered Vanka patches associated with those pressures. Thus, the nonzero graph of \(\Eeul\) identifies the elasticity-induced velocity couplings, and the MAC-grid divergence stencil supplies the geometric closure of the patch. CAV thereby follows the organizational principle of vertex-star relaxation---forming larger patches as unions of simpler geometric patches---but determines these unions using coupling information encoded in the nonzero pattern of \(\Eeul\) rather than finite-element mesh incidence. The multigrid smoother studied here applies multiplicative Schwarz relaxation to the resulting patch family. Because the spreading, interpolation, and Lagrangian elasticity matrices are local in the sense that their stencil sizes remain bounded under mesh refinement, CAV patch sizes are likewise bounded independently of mesh resolution. Consequently, under grid refinement, each multigrid cycle has linear complexity in the total number of Eulerian degrees of freedom.

The numerical study examines the empirical spectra of the implemented multiplicative cycle, together with the performance of CAV-preconditioned FGMRES for target-point, membrane, and beam force laws. The preconditioned iteration counts are nearly independent of grid resolution and remain modest across the tested ranges of stiffness and Reynolds number. A nonlinear benchmark inspired by the flexible-fiber drag-reduction experiments and theory of Alben, Shelley, and Zhang~\cite{alben2002} tests the preconditioner in a finite-Reynolds-number channel with inflow and traction boundary conditions. Under simultaneous Eulerian and Lagrangian refinement, the average number of FGMRES iterations per Newton solve increases modestly, from $8.65$ to $9.53$, while the maximum CAV patch size remains bounded. Appendix~\ref{sec:appendix} offers a heuristic interpretation of these results using an additive Schwarz model of the momentum block alone. This model relates CAV patch selection to a local representation of the elastic near-kernel and suggests a natural additive extension, but it does not constitute a convergence theory for the pressure-coupled multiplicative multigrid method studied here. To our knowledge, CAV is the first preconditioner for the implicit IB equations to combine linear per-cycle complexity with nearly grid-independent convergence across target-point, membrane, and beam force laws and broad ranges of stiffness and Reynolds number.

\section{The Eulerian saddle-point system}
\label{sec:background}
\label{sec:eulerian_schur}

We discretize velocity and pressure using a uniform \(N\times N\) two-dimensional MAC grid with mesh width \(h=1/N\)~\cite{harlow1965,lebedev1964}. The indices \(i,j=0,\ldots,N-1\) enumerate the unique periodic degrees of freedom. The pressure \(p_{i,j}\) is stored at the cell center \(\xb_{i,j}=\left((i+\tfrac12)h,(j+\tfrac12)h\right)\). The horizontal and vertical velocity components \(u_{i,j}\) and \(v_{i,j}\) are stored at the face locations \(\xb_{i-\frac12,j}=\left(ih,(j+\tfrac12)h\right)\) and \(\xb_{i,j-\frac12}=\left((i+\tfrac12)h,jh\right)\), respectively. Let \(\Uh\) denote the discrete periodic MAC velocity space and let \(\Ph=\{p=(p_{i,j}):\sum_{i,j}p_{i,j}=0\}\) denote the corresponding mean-zero pressure space. Unless stated otherwise, the velocity--pressure operators below act on \(\Uh\times\Ph\). Recall that \(\Db\), \(\Gb\), and \(\Lb\) are the discrete divergence, gradient, and vector Laplacian matrices, respectively. The standard MAC discretization satisfies \(\Db=-\Gb\tran\). The flexible fiber model studied in Subsection~\ref{subsec:asz_fiber} uses the same staggered-grid arrangement, with periodic lateral boundaries, an inflow condition on the top boundary, and a traction condition on the bottom boundary. The traction and inflow boundary conditions are implemented for the MAC grid using the ghost-cell filling routines described by Griffith~\cite{griffith2009}.

We represent the immersed structure by \(M\) Lagrangian markers \(\Xb_k\). We construct \(\Sb[\Xb]\) and \(\Jb[\Xb]\) from the tensor-product regularized delta function \(\delta_h\). For \(\xb=(x,y)\), \(\delta_h(\xb)=h^{-2}\varphi(x/h)\varphi(y/h)\), and \(\varphi\) is the one-dimensional kernel. The force density at marker \(k\) is \(\Fb_k=(F_k^x,F_k^y)\), and its Lagrangian quadrature weight is \(\omega_k\). The Eulerian force density \(\fb=\Sb[\Xb]\Fb\) and Lagrangian marker velocity \(\Ub=\Jb[\Xb]\ub\) are given by
\begin{align*}
\fb_{i,j}
&=\sum_{k=1}^{M}
\begin{bmatrix}
F_k^x\,\delta_h\!\left(\xb_{i-\frac12,j}-\Xb_k\right)\\
F_k^y\,\delta_h\!\left(\xb_{i,j-\frac12}-\Xb_k\right)
\end{bmatrix}
\omega_k,\\
\Ub_k
&=\sum_{i=0}^{N-1}\sum_{j=0}^{N-1}
\begin{bmatrix}
u_{i,j}\,\delta_h\!\left(\xb_{i-\frac12,j}-\Xb_k\right)\\
v_{i,j}\,\delta_h\!\left(\xb_{i,j-\frac12}-\Xb_k\right)
\end{bmatrix}h^2.
\end{align*}
With the quadrature-weighted discrete Eulerian and Lagrangian inner products \(\langle\cdot,\cdot\rangle_{\mathrm E}\) and \(\langle\cdot,\cdot\rangle_{\mathrm L}\), these definitions satisfy \(\langle\ub,\Sb[\Xb]\Fb\rangle_{\mathrm E}=\allowbreak\langle\Jb[\Xb]\ub,\Fb\rangle_{\mathrm L}\) for every grid velocity \(\ub\) and marker force density \(\Fb\). Hence, \(\Jb[\Xb]=\Sb[\Xb]\adj\). Throughout this study, we use Peskin's four-point kernel~\cite{peskin2002}, so each Lagrangian marker couples to a \(4\times4\) stencil of grid faces for each velocity component.

For the semi-implicit backward Euler discretization used throughout, we evaluate both matrices at the known configuration \(\Xb^n\): \(\Sb^n=\Sb[\Xb^n]\) and \(\Jb^n=\Jb[\Xb^n]\). The scheme is
\begin{align}
\rho \frac{\ub^{n+1} - \ub^n}{\dt} &= -\rho\,\mathcal N[\ub^n]+\mu\Lb\ub^{n+1} - \Gb p^{n+1} + \Sb^n\Fb^{n+1},\label{eq:im_disc_mom}\\
\Db\ub^{n+1} &= \zerob,\notag\\
\frac{\Xb^{n+1} - \Xb^n}{\dt} &= \Jb^n\ub^{n+1}.\label{eq:im_disc_advec}
\end{align}
The nonlinear convective operator \(\mathcal N[\ub]\) approximates \((\ub\cdot\nabla)\ub\). We approximate all Eulerian spatial derivatives by standard second-order centered differences. Because the two velocity components are stored on different sets of cell faces, evaluating \(\mathcal N[\ub]\) requires interpolating each component to the locations of the other. Explicit treatment of convection imposes the convective CFL restriction \(\dt=\mathcal O(h)\). Implicit treatment of the immersed boundary force removes the generally more severe stiffness-dependent stability restriction.

For a linear force law, \(\Fb^{n+1}=\Eb\Xb^{n+1}\), in which \(\Eb\) is the discrete Lagrangian elasticity matrix. For example, using the linear spring-tension law \(\Fb=\kappa_{\mathrm s}\,\partial^2\Xb/\partial s^2\) on a closed membrane with uniform Lagrangian spacing \(\ds\) gives \begin{equation}
\label{eq:disc_force_law}
\Fb_k^{n+1} = \frac{\kappa_{\mathrm s}}{\ds^2}\left(\Xb^{n+1}_{k+1} + \Xb^{n+1}_{k-1} - 2\Xb^{n+1}_{k}\right).
\end{equation}
For a nonlinear discrete force law \(\Fb=\mathcal F[\Xb]\), the Lagrangian elasticity matrix is the Jacobian \(\Eb[\Xb]=D\mathcal F[\Xb]\). At each Newton iteration, we evaluate \(\Eb[\Xb]\) at the current iterate, and the linearized force increment satisfies \(\delta\Fb=\Eb[\Xb]\,\delta\Xb\).

Equations~\eqref{eq:im_disc_mom}--\eqref{eq:im_disc_advec} give the extended Eulerian--Lagrangian saddle-point system for a linear force law:
\begin{align}
    \label{eq:ext_saddle_sys}
    \underbrace{\begin{bmatrix}
    \Azero & \Gb & -\Sb^n\Eb \\
    -\Db & \zeromat & \zeromat \\
    -\Jb^n & \zeromat & \dt^{-1}\Ib
    \end{bmatrix}}_{\Kb^n}
    \begin{bmatrix}
        \ub^{n+1} \\
        p^{n+1} \\
        \Xb^{n+1}
    \end{bmatrix}
    =
    \begin{bmatrix}
    \bb \\
    \zerob \\
    \dt^{-1}\Xb^n
    \end{bmatrix},
\end{align}
in which \(\bb=(\rho/\dt)\ub^n-\rho\mathcal N[\ub^n]\).

Eliminating the Lagrangian variables from Equation~\eqref{eq:ext_saddle_sys} gives the reduced Eulerian system for velocity and pressure:
\begin{align}
    \label{eq:Eul_schur}
    \underbrace{\begin{bmatrix}
        \Ab^n & \Gb \\
        -\Db & \zeromat
    \end{bmatrix}}_{\KE^n}
        \begin{bmatrix}
            \ub^{n+1} \\
            p^{n+1}
    \end{bmatrix}
    =
    \begin{bmatrix}
        \bb + \Sb^n\Eb\Xb^n \\
        \zerob
    \end{bmatrix}.
\end{align}
Here, \(\Eeul^n=\Sb^n\Eb\Jb^n\) and \(\Ab^n=\Azero-\dt\,\Eeul^n\). Recall that \(\Azero=(\rho/\dt)\Ib-\mu\Lb\), so both \(\Azero\) and \(\Ab^n\) depend on the time-step size. If the time level is immaterial, we omit the superscript \(n\) from \(\Eeul\), \(\Ab\), and \(\KE\).

Eliminating the Eulerian variables from Equation~\eqref{eq:ext_saddle_sys} instead produces a Lagrangian Schur complement containing the inverse of the discrete Stokes operator~\cite{newren2007a,newren2007}. This Schur complement is generally dense and has proved difficult to precondition effectively, particularly in the time-dependent setting. Further eliminating velocity or pressure from Equation~\eqref{eq:Eul_schur} produces either a dense projected velocity operator or a pressure Schur complement containing \(\Ab^{-1}\)~\cite{newren2007a}. These reductions destroy the sparse saddle-point structure of \(\KE\). We therefore precondition \(\KE\) directly, exploiting its sparse Cartesian-grid structure with geometric multigrid.

\section{Standard and generalized Vanka patches}
\label{sec:cav}

\subsection{Vanka and box-relaxation patches}

Vanka~\cite{vanka1986} introduced pressure-centered patch relaxation for the incompressible Navier--Stokes equations. All cell indices are interpreted periodically modulo \(N\). The global set of velocity degrees of freedom is \(\mathcal U\). The velocity degrees of freedom on the boundary of cell \((i,j)\) and the corresponding standard pressure-centered Vanka patch are
\begin{equation*}
\mathcal U_{i,j}
=
\{u_{i,j},u_{i+1,j},v_{i,j},v_{i,j+1}\},
\qquad
\mathcal V_{i,j}
=
\{p_{i,j}\}\cup\mathcal U_{i,j}.
\end{equation*}
In a Vanka smoothing sweep, solving the \(5\times5\) local saddle-point system on \(\mathcal V_{i,j}\) updates the pressure unknown and the four velocity unknowns in its discrete divergence equation simultaneously.

Larger generalized Vanka patches can be formed as unions of standard patches. For any set \(C\) of cell indices, the union of their standard Vanka patches is
\begin{equation*}
\mathcal V(C)
=
\bigcup_{(i,j)\in C}\mathcal V_{i,j}.
\end{equation*}
Box relaxation constructs larger patches by choosing \(C\) geometrically.
The box construction \(\mathsf B(b,o)\) partitions the MAC grid into nonoverlapping \(b\times b\) rectangular cell blocks and expands each block by \(o\) cell layers in every coordinate direction~\cite{guy2015}. For a block whose lower-left cell is \((i,j)\), its expanded rectangular cell set and individual box patch are
\begin{equation*}
\mathcal R_{i,j}^{b,o}
=
\{(i+r,j+s):-o\le r,s\le b+o-1\},
\qquad
\mathcal V_{i,j}^{b,o}
=
\mathcal V(\mathcal R_{i,j}^{b,o}).
\end{equation*}
As \((i,j)\) ranges over the lower-left cells of the blocks, \(\mathsf B(b,o)\) produces the family of patches \(\mathcal V_{i,j}^{b,o}\). Thus, \(\mathsf B(1,0)\) produces the standard Vanka family \(\{\mathcal V_{i,j}\}\). All box comparisons below use \(\mathsf B(4,2)\).

\subsection{Coupling-aware Vanka patches}
\label{subsec:cav_patches}

CAV provides an alternative generalization of standard Vanka patches that uses the nonzero pattern of \(\Eeul\) to determine which patches are combined. CAV forms an expanded velocity set consisting of the velocities in a standard pressure-centered Vanka patch together with every velocity degree of freedom coupled to one of those velocities through a nonzero entry of \(\Eeul\). If the expansion does not enlarge the standard velocity set, CAV retains the original standard patch. Otherwise, it takes the union of the standard Vanka patches associated with the cells incident to the expanded velocity set. To state this construction precisely, the cell-incidence map \(\mathcal C\) associates each velocity degree of freedom with its two incident cells:
\begin{align*}
\mathcal C(u_{i,j})
&=\{(i-1,j),(i,j)\},&
\mathcal C(v_{i,j})
&=\{(i,j-1),(i,j)\}.
\end{align*}
For any \(\mathcal W\subseteq\mathcal U\), this map extends by set union:
\begin{equation*}
\mathcal C(\mathcal W)
=
\bigcup_{w\in\mathcal W}\mathcal C(w).
\end{equation*}
Equivalently, \(\mathcal C(\mathcal W)\) is the set of cells whose corresponding rows of \(\Db\) involve at least one velocity degree of freedom in \(\mathcal W\).

\begin{figure}[t]
    \centering
    \includegraphics[width=\textwidth]{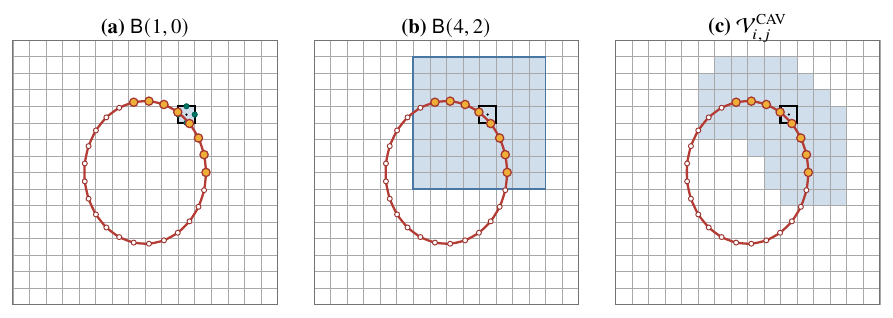}
    \caption{Representative \(\mathsf B(1,0)\), \(\mathsf B(4,2)\), and CAV patches generated from the same pressure seed, marked by the black square, for a linearly elastic membrane on an \(N=16\) MAC grid. Blue shading identifies the cells whose standard Vanka patches are united, and red circles denote the Lagrangian markers discretizing the membrane. Gold markers identify the Lagrangian nodes that mediate the elasticity-induced velocity couplings associated with the seed patch. The same nodes are highlighted in all three panels to provide a common reference, although only CAV uses this coupling information to determine the patch. Because the four-point regularized delta kernel couples each marker to Eulerian velocity degrees of freedom across four grid cells in each coordinate direction, the CAV patch extends beyond the cells containing the highlighted markers.}
    \label{fig:patch_constructions}
\end{figure}

For each cell \((i,j)\), the expanded velocity set is
\begin{equation*}
\widehat{\mathcal U}_{i,j}
=
\mathcal U_{i,j}
\cup
\left\{
k\in\mathcal U:
(\Eeul)_{k\ell}\ne0
\ \text{or}\
(\Eeul)_{\ell k}\ne0
\text{ for some }\ell\in\mathcal U_{i,j}
\right\}.
\end{equation*}
The CAV patch associated with cell \((i,j)\) is
\begin{equation}
\label{eq:cav_patch_set}
\mathcal V_{i,j}^{\text{CAV}}
=
\begin{cases}
\mathcal V_{i,j},
&
\widehat{\mathcal U}_{i,j}=\mathcal U_{i,j},
\\[2mm]
\mathcal V\!\left(\mathcal C(\widehat{\mathcal U}_{i,j})\right),
&
\text{otherwise}.
\end{cases}
\end{equation}
Equation~\eqref{eq:cav_patch_set} associates one CAV patch with every cell. The first branch retains the standard Vanka patch if the expanded velocity set is unchanged. Otherwise, \(\mathcal C\) identifies the cells incident to the expanded velocity set, and \(\mathcal V\) forms the union of their standard Vanka patches. \Cref{fig:patch_constructions} compares the three patch constructions for a linearly elastic membrane.

\subsection{Multiplicative patch relaxation}
\label{subsec:patch_relaxation}

Given ordered patches \(\mathcal V_1,\ldots,\mathcal V_m\), a standard multiplicative smoother uses the matrix \(\Rb_p\) to restrict a global velocity--pressure vector to \(\mathcal V_p\), and \(\Kb_{\text{E},p}=\Rb_p\KE\Rb_p\tran\) is the corresponding local matrix. Each sweep begins with \(\rb\leftarrow\bb-\KE\wb\). For \(p=1,\ldots,m\), the smoother applies
\begin{align*}
\deltab_p&\leftarrow\Kb_{\text{E},p}^{-1}\Rb_p\rb,\\
\wb&\leftarrow\wb+\Rb_p\tran\deltab_p,\\
\rb&\leftarrow\rb-\KE\Rb_p\tran\deltab_p.
\end{align*}
After the sweep, we remove the mean pressure. Our CAV implementation associates one patch with each Cartesian grid cell and visits the patches in lexicographic order. For the box comparisons, the geometric \(\mathsf B(b,o)\) family replaces the CAV patch family, and we visit its patches in lexicographic order of the lower-left cells of their unexpanded blocks. Our numerical tests use the box-relaxation smoother of Guy et~al.~\cite{guy2015} with restricted injection. We solve the local system associated with each box patch \(\mathcal V_{i,j}^{4,2}\) and inject the resulting correction only into the degrees of freedom associated with its unexpanded \(4\times4\) rectangular cell block. These cell blocks form a nonoverlapping partition, whereas adjacent blocks share their face-centered velocity degrees of freedom. Except for this restricted injection, both constructions use the same relaxation algorithm.

For both patch families, we apply one pass of max-norm row and column equilibration to each local saddle-point matrix before computing its LU factorization to mitigate poor conditioning in finite-precision arithmetic. We scale the residual and correction consistently, so this equilibration leaves the local Schwarz correction unchanged in exact arithmetic.

The cost of patch relaxation depends on the local dimensions \(n_p=|\mathcal V_p|\). If each local matrix is factorized independently, the setup cost is \(\mathcal O(\sum_p n_p^3)\), and one sweep with the stored factors costs \(\mathcal O(\sum_p n_p^2)\). Equilibration adds \(\mathcal O(\sum_p n_p^2)\) work during setup and \(\mathcal O(\sum_p n_p)\) work per sweep, so it does not change these bounds. For CAV, the fixed support of the regularized delta function, the fixed-width Lagrangian force stencils, and our choice \(\Delta s=\mathcal O(h)\) keep \(n_p\) bounded under refinement. For \(\mathsf B(b,o)\) with fixed \(b\) and \(o\), the local dimensions are fixed. Each smoothed grid level contains one CAV patch per cell and one box patch per \(b\times b\) rectangular cell block, so both patch counts scale linearly with the number of cells.

\section{Geometric multigrid preconditioner}
\label{sec:multigrid}

We define the geometric multigrid preconditioner on a grid hierarchy formed by halving the number of MAC cells in each coordinate direction. Between adjacent levels \(\ell-1\) and \(\ell\), the velocity prolongation \(\Pb_{\ell-1}^{\ell}\) is defined by lowest-order Raviart--Thomas interpolation on quadrilaterals~\cite{raviart1977}, and the velocity restriction is its discrete \(L^2\)-adjoint, \(\Rb_{\ell}^{\ell-1}=(\Pb_{\ell-1}^{\ell})\adj=\tfrac14(\Pb_{\ell-1}^{\ell})\tran\). This scaling accounts for the refinement ratio in two spatial dimensions and preserves constant velocities under restriction. We prolong pressure by bilinear interpolation and restrict it by cell averaging.

On each level, we rediscretize the standard time-dependent Stokes part of the operator: the discrete divergence and gradient matrices \(\Db\) and \(\Gb\), together with the mass and vector Laplacian terms that form \(\Azero\). In contrast, we assemble the Eulerian elasticity matrix on the finest level \(\ell_{\text{max}}\) and use Galerkin projection to construct the corresponding matrices on coarser levels. The matrix \(\Eeul^{(\ell)}\) is the Eulerian elasticity matrix on level \(\ell\). Starting from \(\Eeul^{(\ell_{\text{max}})}\), the coarse-grid matrices satisfy
\begin{equation*}
\Eeul^{(\ell-1)}
=
\Rb_{\ell}^{\ell-1}
\Eeul^{(\ell)}
\Pb_{\ell-1}^{\ell},
\qquad \ell=\ell_{\text{max}},\ldots,1.
\end{equation*}
On every level except the coarsest, we construct the selected patch family and compute the corresponding local LU factorizations. For CAV, the patch construction uses the sparsity pattern of the Galerkin-projected Eulerian elasticity matrix on that level. We solve the coarsest saddle-point system using a precomputed direct factorization. Each V-cycle applies one pre-sweep and one post-sweep. The per-level cost bounds from Subsection~\ref{subsec:patch_relaxation} and the geometric decrease in cell counts between levels imply that the patch-factorization setup and relaxation work over the full hierarchy both scale linearly with the number of fine-grid Eulerian unknowns.

\section{Numerical results}
\label{sec:results}

We conduct three sets of numerical experiments. The linearly elastic membrane experiments compare CAV with the $\mathsf B(1,0)$ and $\mathsf B(4,2)$ patch constructions and examine two-grid spectra and FGMRES convergence with multigrid V-cycle preconditioning under mesh refinement~\cite{guy2015}. The bending-resistant closed-beam and target-point experiments repeat the refinement study for force laws that produce different Eulerian coupling patterns. The nonlinear flow past a flexible fiber benchmark tests CAV preconditioning of Newton Jacobian systems under simultaneous Eulerian and Lagrangian refinement.

Unless stated otherwise, the linear experiments use a zero initial guess and right-preconditioned FGMRES. Residuals are measured in the Euclidean norm, and convergence is declared if the relative residual satisfies \(\|\rb_k\|_2/\|\rb_0\|_2\le10^{-10}\). The membrane solves are limited to a maximum of \(60\) iterations, whereas the beam and target-point solves are limited to \(150\). We do not use any restarts. The nonlinear solver tolerances are specified in \Cref{subsec:asz_fiber}.

\subsection{Linear elastic membrane forces}
\label{subsec:membrane}

The membrane experiments use the three-point discretization of the linear spring-tension law in Equation~\eqref{eq:disc_force_law} on an ellipse centered at \((1/2,1/2)\) with horizontal and vertical semi-axes \(0.23\) and \(0.27\), respectively. At every resolution, the Lagrangian marker spacing is approximately \(h/2\), and the time-step size is \(\Delta t=h/2\). We form each system from the first backward Euler time step with zero initial fluid velocity.

The heuristic momentum-block model in Appendix~\ref{sec:appendix} suggests that CAV patch selection should prevent stiffness-dependent eigenvalues from approaching zero. We compare this qualitative prediction with the implemented pressure-coupled CAV cycle. The operator \(\Bb^{\text{TG}}\) applies one two-grid cycle with one multiplicative pre-sweep, one multiplicative post-sweep, and an exact coarse-grid solve. Its iteration matrix is \(\Tb=\Ib-\Bb^{\text{TG}}\KE\). The spectral radius of \(\Tb\) is \(\rho(\Tb)\). The eigenvalues of the preconditioned matrix \(\Bb^{\text{TG}}\KE\) are \(1-\lambda(\Tb)\). Consequently, clustering of \(\lambda(\Tb)\) near the origin corresponds to clustering of the preconditioned spectrum near one. If these matrices are instead assembled on the full pressure coefficient space, the periodic saddle-point matrix has a constant-pressure null mode, which produces a unit eigenvalue of the full-space iteration matrix and a zero eigenvalue of the full-space preconditioned matrix. Omitting these gauge modes, as we do in the reported spectra, is equivalent to restricting the operators to \(\Uh\times\Ph\).

\begin{table}[t]
\centering
\small
\caption{Spectrum of the implemented pressure-coupled CAV two-grid cycle for the $N=16$ linear elastic membrane test with an $8\times8$ coarse grid, $\mu=10^{-2}$, $\rho=1$, one multiplicative pre-sweep, and one multiplicative post-sweep. The constant-pressure mode is omitted.}
\label{tab:implemented_cav_stiffness_spectrum}
\begin{tabular}{lccc}
\toprule
$\kappa_{\mathrm s}$ & $\rho(\Tb)$
& $\min_{\lambda}|\lambda(\Bb^{\text{TG}}\KE)|$
& $\max_{\lambda}|\lambda(\Bb^{\text{TG}}\KE)|$ \\
\midrule
$10^{2}$ & $0.0816$ & $0.9357$ & $1.0486$ \\
$10^{3}$ & $0.1233$ & $0.9471$ & $1.1233$ \\
$10^{4}$ & $0.3518$ & $0.8470$ & $1.3518$ \\
$10^{5}$ & $0.4684$ & $0.8090$ & $1.4684$ \\
$10^{6}$ & $0.5916$ & $0.8326$ & $1.5916$ \\
\bottomrule
\end{tabular}
\end{table}

We use the spectra in \cref{tab:implemented_cav_stiffness_spectrum} as a qualitative diagnostic for preconditioned FGMRES. As $\kappa_{\mathrm s}$ increases from $10^2$ to $10^6$, the spectrum broadens gradually, revealing a weak stiffness dependence that is not captured by the heuristic additive momentum-block model in Appendix~\ref{sec:appendix}. This difference is not unexpected because the heuristic momentum-block model omits pressure coupling and the multiplicative multilevel construction. In addition, the implemented iteration matrix is generally nonnormal, so its eigenvalues alone do not determine FGMRES convergence. Nevertheless, we emphasize that the observed dependence is mild: over four decades in $\kappa_{\mathrm s}$, every nonzero eigenvalue of the preconditioned matrix has modulus at least $0.8090$, and the largest modulus is $1.5916$. The spectral results therefore support the qualitative prediction that CAV avoids stiffness-dependent near-zero modes. The following experiments provide more direct convergence evidence from FGMRES.

\begin{figure}[t]
\centering
\begin{tabular}{@{}c@{\hspace{0.025\textwidth}}c@{\hspace{0.025\textwidth}}c@{}}
\subcaptionbox{}[0.29\textwidth]{\includegraphics[width=\linewidth]{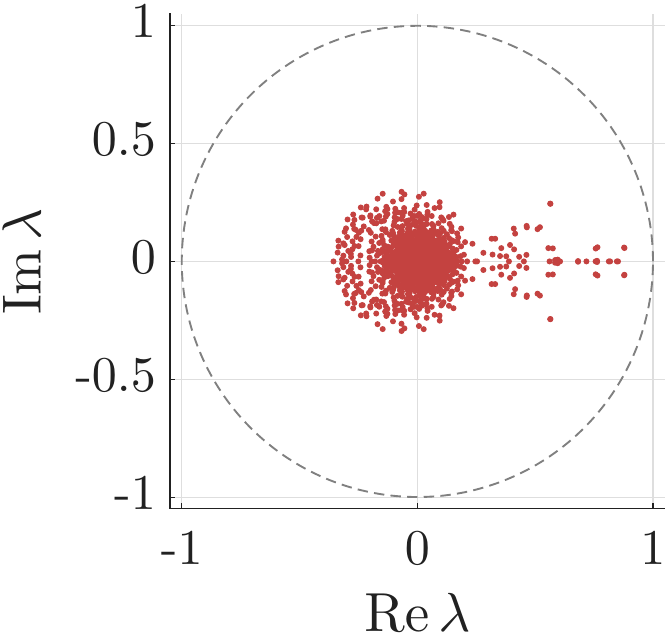}} &
\smash{\raisebox{-0.5\height}{\subcaptionbox{}[0.333\textwidth]{\includegraphics[width=\linewidth]{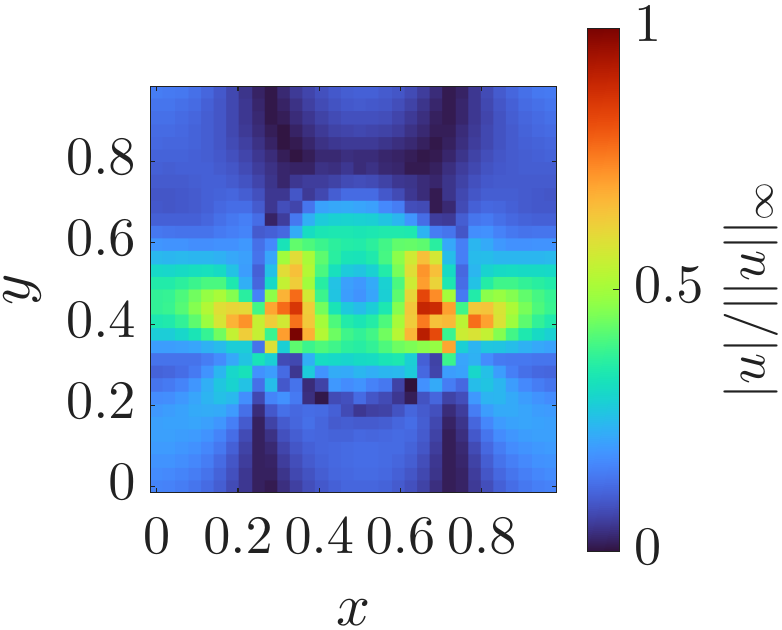}}}} &
\subcaptionbox{}[0.29\textwidth]{\includegraphics[width=\linewidth]{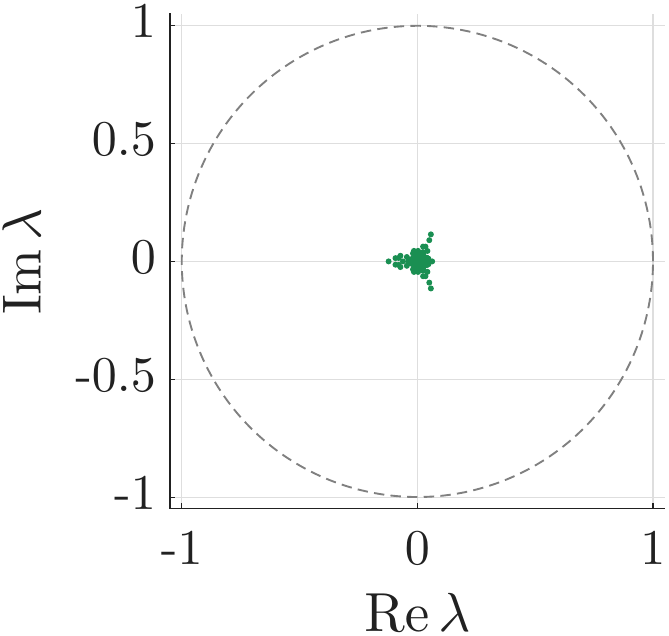}} \\[0.8em]
\subcaptionbox{}[0.29\textwidth]{\includegraphics[width=\linewidth]{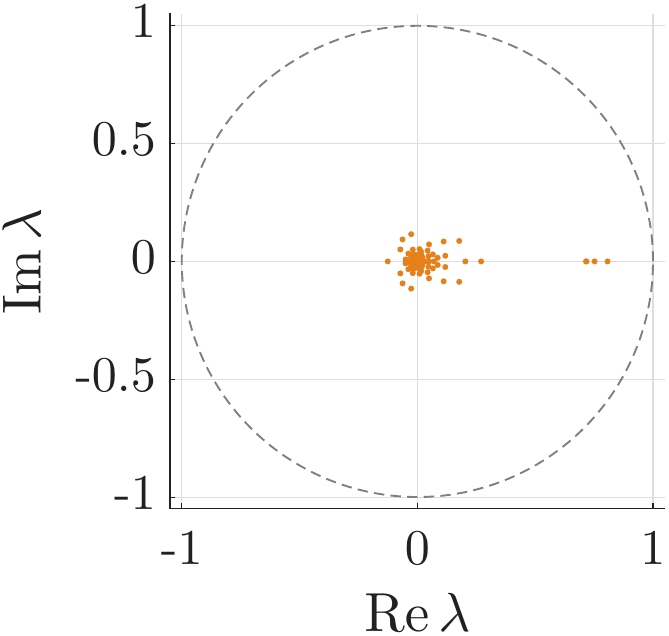}} &&
\subcaptionbox{}[0.29\textwidth]{\includegraphics[width=\linewidth]{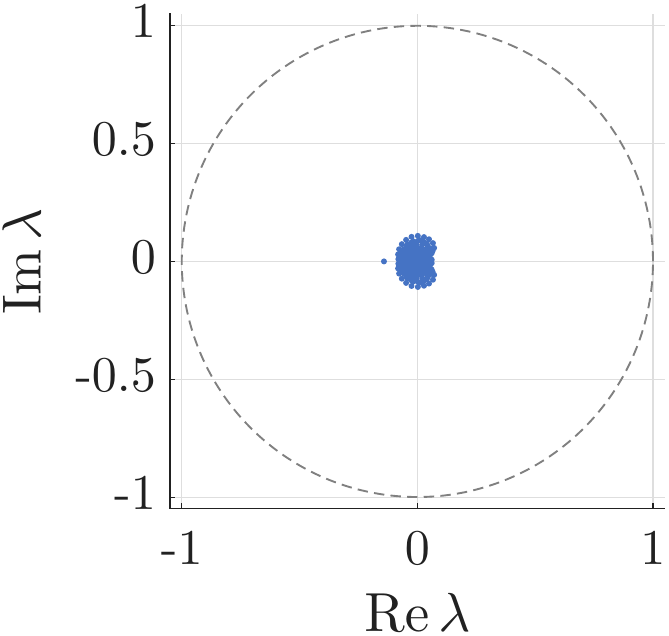}}
\end{tabular}
\caption{Spectra of the two-grid iteration matrices for the linear elastic membrane test ($N = 32$, $16 \times 16$ coarse grid, $\kappa_{\mathrm s} = 10^4$). The constant-pressure mode is omitted, and the dashed curve is the unit circle. (a)~$\mathsf B(1,0)$ relaxation ($\mu = 1$), spectral radius $0.88$. (b)~Magnitude of the horizontal velocity component of the largest-eigenvalue mode from~(a), concentrated near the immersed boundary. (c)~$\mathsf B(4,2)$ relaxation ($\mu = 1$), spectral radius $0.13$. (d)~$\mathsf B(4,2)$ relaxation with $\mu = 10^{-2}$, spectral radius $0.81$. (e)~CAV relaxation with $\mu = 10^{-2}$, spectral radius $0.14$.}
\label{fig:eigenvalue_comparison}
\end{figure}

We next compare CAV with $\mathsf B(1,0)$ and $\mathsf B(4,2)$ patch relaxation for the linear elastic membrane test on the $N=32$ grid, using a $16\times16$ coarse grid. For each relaxation, we assemble the two-grid iteration matrix \(\Tb\) and compute its spectrum. \Cref{fig:eigenvalue_comparison} shows the resulting spectra.

For $\rho = 1$, $\mu = 1$, and $\kappa_{\mathrm s} = 10^4$, $\mathsf B(1,0)$ relaxation yields eigenvalues confined to the unit disk, but they are far from clustered: the spectral radius is $0.88$ (\cref{fig:eigenvalue_comparison}a), and the slowest mode's horizontal velocity component concentrates near the immersed boundary (\cref{fig:eigenvalue_comparison}b), indicating that the standard pressure-centered Vanka patch does not resolve the coupling introduced by the Eulerian elasticity operator. The $\mathsf B(4,2)$ construction performs substantially better in this regime, with spectral radius $0.13$ (\cref{fig:eigenvalue_comparison}c). Reducing the viscosity to $\mu = 10^{-2}$, however, degrades this box construction sharply, with the spectral radius growing to $0.81$ (\cref{fig:eigenvalue_comparison}d). The corresponding CAV iteration matrix in the same low-viscosity, high-stiffness regime keeps the eigenvalues tightly clustered near the origin, with spectral radius $0.14$ (\cref{fig:eigenvalue_comparison}e).

We then assess preconditioned FGMRES convergence under grid refinement. Using the same membrane configuration, force law, Lagrangian spacing, and time-step size, we compare $\mathsf B(4,2)$ and CAV patch relaxation as multigrid V-cycle right preconditioners for FGMRES with $\kappa_{\mathrm s} = 10^4$, $\mu = 10^{-2}$, and $\rho = 1$ on the $N=16$, $32$, and $64$ grids. Each V-cycle uses the $8\times8$ grid as its coarsest level.

\begin{figure}[t]
    \centering
    \includegraphics[width=0.75\textwidth]{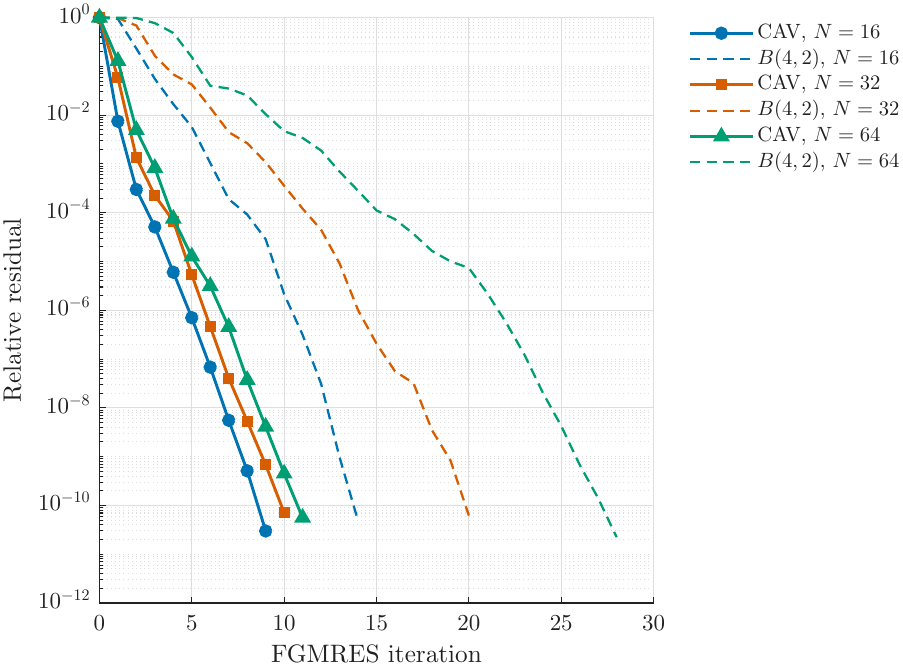}
    \caption{Convergence of FGMRES with multigrid V-cycle preconditioning for the linear elastic membrane test ($\kappa_{\mathrm s} = 10^4$, $\mu = 10^{-2}$, $N = 16$, $32$, and $64$). CAV patch relaxation is shown with solid, filled markers, and $\mathsf B(4,2)$ patch relaxation is shown with dashed curves.}
    \label{fig:gmres_convergence}
\end{figure}

CAV reduces the relative residual by ten orders of magnitude in $9$, $10$, and $11$ iterations at $N = 16$, $32$, and $64$ (\cref{fig:gmres_convergence}). The $\mathsf B(4,2)$ preconditioner converges but deteriorates under grid refinement, requiring $14$, $20$, and $28$ iterations, respectively.

The nearly constant CAV iteration counts demonstrate grid-scalable preconditioning for the membrane test. By contrast, the $\mathsf B(4,2)$ iteration count doubles from $14$ to $28$ over the same refinement sequence, showing that enlarging Vanka patches by a fixed amount does not yield a grid-scalable preconditioner for this time-dependent membrane model. Together with the two-grid spectra, these results indicate that CAV captures the elasticity-induced couplings that govern convergence in the low-viscosity, high-stiffness regime.

\subsection{Beam and target-point forces}
\label{subsec:other_forcings}

The stiffness of the discrete beam operator grows rapidly under grid refinement, making the beam test a particularly demanding assessment of preconditioner robustness. We compare the beam and target-point tests to evaluate the same CAV construction for two force laws with sharply different coupling patterns.

The continuum Euler--Bernoulli bending force density is
\begin{equation*}
\Fb=-\kappa_{\mathrm b}\frac{\partial^4\Xb}{\partial s^4}.
\end{equation*}
Target-point tethering instead applies the markerwise force
\begin{equation*}
\Fb_k(t)=\kappa_{\mathrm t}\left(\Xb_k^{\mathrm{target}}(t)-\Xb_k(t)\right).
\end{equation*}
The prescribed target position \(\Xb_k^{\mathrm{target}}(t)\) varies in time.

For the beam tests, the immersed structure is the closed curve with polar radius \(r(\theta)=0.23+0.035\cos(3\theta)\), centered in the periodic unit square. Its Lagrangian markers are spaced approximately one Eulerian mesh width apart. We discretize the bending law using a standard centered five-point difference stencil, which couples each marker to two neighbors on either side. The stencil coefficients are proportional to \(\kappa_{\mathrm b}\Delta s^{-4}\); because we choose \(\Delta s\) proportional to \(h\), the discrete beam stiffness grows proportionally to \(h^{-4}\) under refinement. We use \(\kappa_{\mathrm b}=1\) and \(\mu=10^{-2}\), retaining the low-viscosity regime used in the membrane comparison.

The target-point test uses two rows of markers initially coincident with their prescribed targets at \(y=0.25\) and \(y=0.75\). The lower and upper targets move horizontally with velocities \(-U\) and \(+U\), respectively, with \(U=0.05\). The target-point elasticity matrix is \(-\kappa_{\mathrm t}\Ib\) and is therefore diagonal in the Lagrangian marker index, introducing no intermarker coupling. We set \(\kappa_{\mathrm t}=10^6\) and \(\mu=1\). Because the purpose of this experiment is to assess solver scalability rather than convergence to an exactly constrained moving wall, we hold \(\kappa_{\mathrm t}\) fixed under grid refinement. Thus, every resolution uses the same strongly penalized target-point force law.

In both experiments, \(\rho=1\) and \(\Delta t=h/2\), and we solve the system arising from the first backward Euler time step with zero initial fluid velocity. We compare CAV, $\mathsf B(4,2)$ patch relaxation, and unpreconditioned FGMRES on the $N=32$, $64$, and $128$ grids. FGMRES stops when the relative residual reaches $10^{-10}$, and each V-cycle uses an $8\times8$ grid at the coarsest level.

\begin{figure}[t]
\centering
\begin{subfigure}{0.48\textwidth}
    \centering
    \includegraphics[width=\textwidth]{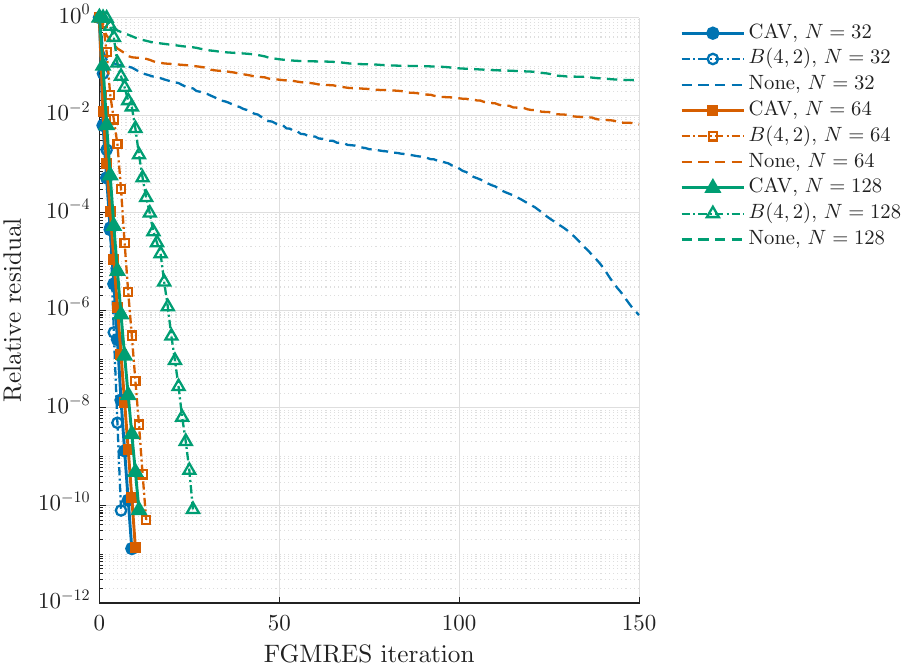}
    \caption{Closed beam, $\kappa_{\mathrm b} = 1$, $\mu = 10^{-2}$.}
\end{subfigure}
\hfill
\begin{subfigure}{0.48\textwidth}
    \centering
    \includegraphics[width=\textwidth]{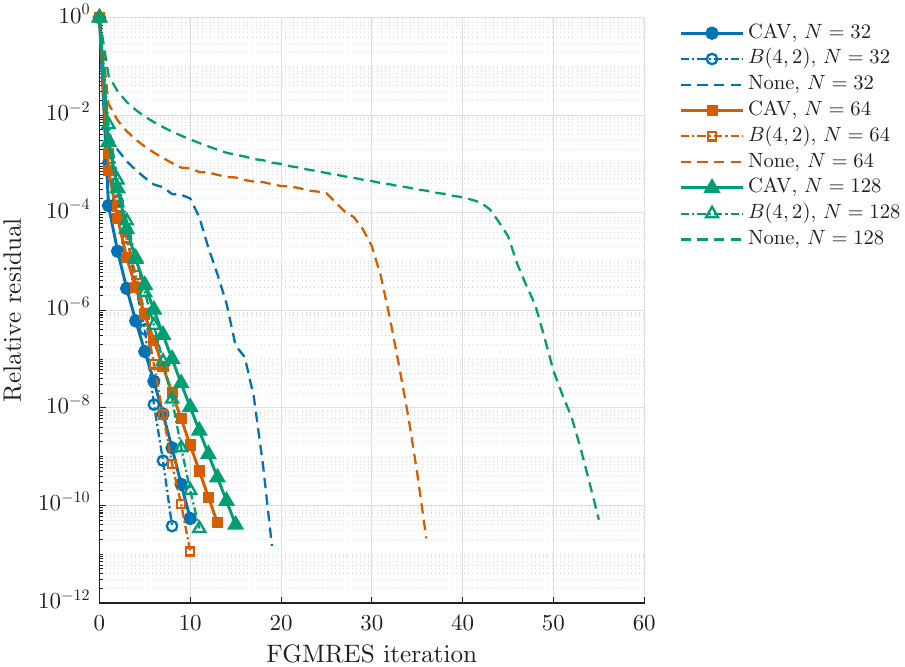}
    \caption{Target-point walls, $\kappa_{\mathrm t} = 10^{6}$, $\mu = 1$.}
\end{subfigure}
\caption{Convergence of FGMRES for the beam and target-point tests defined in Subsection~\ref{subsec:other_forcings} ($N = 32$, $64$, and $128$). The curves show CAV V-cycle preconditioning (solid, filled markers), $\mathsf B(4,2)$ V-cycle preconditioning (dash-dotted, open markers), and no preconditioning (dashed).}
\label{fig:other_forcings_convergence}
\end{figure}

\begin{table}[t]
\centering
\caption{FGMRES iterations required to reduce the relative residual below $10^{-10}$ for the isolated beam and target-point experiments. A dash indicates failure to converge within $150$ iterations.}
\label{tab:other_forcings_iters}
\begin{tabular}{ccccccc}
\toprule
 & \multicolumn{3}{c}{Beam} & \multicolumn{3}{c}{Target points} \\
\cmidrule(lr){2-4}\cmidrule(lr){5-7}
$N$ & CAV & $\mathsf B(4,2)$ & None & CAV & $\mathsf B(4,2)$ & None \\
\midrule
$32$  & $9$  & $6$  & --- & $10$ & $8$  & $19$ \\
$64$  & $10$ & $13$ & --- & $13$ & $10$ & $36$ \\
$128$ & $11$ & $26$ & --- & $15$ & $11$ & $55$ \\
\bottomrule
\end{tabular}
\end{table}

\Cref{fig:other_forcings_convergence} and \cref{tab:other_forcings_iters} summarize the results. For the beam test, unpreconditioned FGMRES does not reach the prescribed tolerance within $150$ iterations at any resolution. The $\mathsf B(4,2)$ preconditioner requires $6$, $13$, and $26$ iterations on the three grids, whereas CAV requires $9$, $10$, and $11$. Thus, the geometric construction gives the lowest count on the coarsest grid but appears to lose robustness under refinement. The rapid growth of the discrete beam stiffness makes the Eulerian couplings that cross fixed geometric subdomain boundaries increasingly important under refinement. CAV identifies these couplings from the assembled operator and retains nearly constant iteration counts. For the target-point test, both V-cycle preconditioners remain effective. The $\mathsf B(4,2)$ counts are $8$, $10$, and $11$, the CAV counts are $10$, $13$, and $15$, and the unpreconditioned counts increase from $19$ to $55$. The target-point operator is diagonal in the Lagrangian index, so a fixed geometric box captures its more localized coupling pattern. These results show the primary benefit of CAV for force laws with extended intermarker coupling, whereas CAV remains effective for the diagonal target-point force law.

We emphasize that the performance of box relaxation depends on the box size and overlap. Enlarging the box allows each local solve to include more of the coupling induced by the beam operator and may reduce the resulting iteration count. However, this improvement comes at the cost of enlarging every local system, including those away from the immersed boundary where the elasticity operator introduces no coupling. Moreover, suitable choices of $b$ and $o$ depend on the force law and discretization. We use $\mathsf B(4,2)$ throughout and do not retune the box construction for individual tests. Thus, $\mathsf B(4,2)$ serves as a fixed geometric reference rather than an optimized box method for each test. The comparison highlights the distinction between a geometric neighborhood chosen in advance and the CAV neighborhood determined from the sparsity pattern of the assembled Eulerian elasticity operator.

\subsection{Flow past a flexible fiber}
\label{subsec:asz_fiber}

We conclude with a nonlinear benchmark inspired by the flexible-fiber drag-reduction experiments and theory of Alben, Shelley, and Zhang~\cite{alben2002}. In their experiments, an initially straight flexible fiber was placed transverse to a descending soap-film flow and clamped at its midpoint. The downward flow loads the fiber and bends its two free ends downstream. Their mathematical model used steady, inviscid, free-streamline theory to describe this high-Reynolds-number experiment. We instead solve the time-dependent, viscous, incompressible Navier--Stokes equations with an immersed boundary discretization. Our goal here is not a quantitative reproduction of the soap-film experiment. Instead, we use it to test the CAV preconditioner in a nonlinear setting, where it preconditions a Newton solver rather than a single linear system. Throughout this section we use the centimeter--gram--second system of units.

The computational domain is $\Omega=(0,2)^2$, periodic in the horizontal direction. The fiber has length $L=0.375$, is initially horizontal, and is centered at $(1,1.5)$. We take $\rho=1$ and $\mu=1.875\times10^{-3}$. A uniform downward velocity is prescribed at the upper boundary, with magnitude increasing smoothly from zero to $U=1$ over $0 \le t \le 0.5625$, i.e., $1.5$ convective time units $L/U$. With this choice of parameters, the flow has a Reynolds number equal to $200$. At the lower boundary, we use the backflow-stabilized traction condition
\begin{equation*}
\bbsigma\bm n=\frac{\rho}{2}\min(\ub\cdot\bm n,0)\ub,
\end{equation*}
in which $\bbsigma=-p\bm I+\mu(\nabla\ub+(\nabla\ub)\tran)$ and $\bm n$ is the outward unit normal. The correction vanishes under outflow and dissipates kinetic energy when the normal velocity reverses~\cite{bruneau1996}.

The fiber mechanics are defined by discrete energy laws. The first is a discrete bending energy that models the fiber's flexibility. The remaining energies enforce the fiber's kinematic constraints: inextensibility, and the midpoint's position and tangent. For the uniform marker spacing $\ds = h/2$ and the straight reference configuration used here, the link joining markers $k$ and $k+1$ is $\Dlink_k=\Xb_{k+1}-\Xb_k$. The discrete bending energy is
\begin{equation*}
\mathcal E_{\mathrm b}[\Xb]
=\frac{EI}{2\ds^5}
\sum_{k=1}^{M-2}
\left(\Dlink_k\times\Dlink_{k+1}\right)^2,
\end{equation*}
in which $\Dlink_k\times\Dlink_{k+1}$ denotes the scalar cross product. It vanishes when adjacent links are parallel, as they are in the straight reference configuration. To enforce the inextensibility constraint, the discrete stretching energy is
\begin{equation*}
\mathcal E_{\mathrm s}[\Xb]
=\frac{K_{\mathrm s}}{2\ds}
\sum_{k=1}^{M-1}
\left(\lVert\Dlink_k\rVert-\ds\right)^2.
\end{equation*}

Let $m$ denote the fiber's midpoint marker, let $\Xb_{\mathrm c}=(1,1.5)$ denote its prescribed position, and let $\bm e_x=(1,0)$. We impose the midpoint clamp of the fiber using the additional penalty energies
\begin{align*}
\mathcal E_{\mathrm p}[\Xb]
&=\frac{K_{\mathrm p}}{2}
\left\lVert\Xb_m-\Xb_{\mathrm c}\right\rVert^2,\\
\mathcal E_{\mathrm t}[\Xb]
&=\frac{K_{\mathrm t}}{2}
\left\lVert
\frac{\Xb_{m+1}-\Xb_{m-1}}{2\ds}-\bm e_x
\right\rVert^2.
\end{align*}
The position penalty keeps the midpoint near $\Xb_{\mathrm c}$, while the tangent penalty keeps the centered tangent horizontal at the clamp. The two fiber ends remain free. The total Lagrangian force density is
\begin{equation*}
\Fb[\Xb]
=-\frac{1}{\ds}\nabla_{\Xb}
\left(
\mathcal E_{\mathrm b}
+\mathcal E_{\mathrm s}
+\mathcal E_{\mathrm p}
+\mathcal E_{\mathrm t}
\right).
\end{equation*}
The factor $\ds^{-1}$ converts the gradient of the discrete energy into a Lagrangian force density because the marker quadrature weight is $\ds$.

The bending stiffness is held fixed at \(EI=2.63672\times10^{-2}\). The stretching, midpoint-position, and midpoint-tangent penalties approximate the exact inextensibility and clamping constraints. To approach the fully constrained formulation as \(h\to0\), the corresponding stiffness parameters must increase under grid refinement. We therefore scale all three parameters in proportion to \(h^{-1}\). At \(N=64\), we set \(K_{\mathrm s}=375\), \(K_{\mathrm p}=1200\), and \(K_{\mathrm t}=7.03125\), and double each value whenever the grid spacing is halved.

Consistent with the semi-implicit backward Euler discretization, the nonlinear fiber force is evaluated at the updated configuration, while spreading and interpolation are evaluated at the beginning of the time step. To avoid cancellation errors in computing the finite differences needed for the bending force, we update the fiber displacement and link vectors incrementally and evaluate the bending force directly from the links. Newton's method produces the Eulerian velocity--pressure system in Equation~\eqref{eq:intro_saddle}, with $\Eb$ replaced by the tangent of the nonlinear fiber force. We rebuild the CAV hierarchy for each Newton linear system and apply one pre-sweep and one post-sweep on every level. The coarsest grid contains $8\times8$ cells. FGMRES terminates when the relative residual is below $10^{-9}$. Newton's method terminates when either the relative residual or absolute residual is below $10^{-8}$.

\begin{figure}[!t]
\centering
\includegraphics[width=\textwidth]{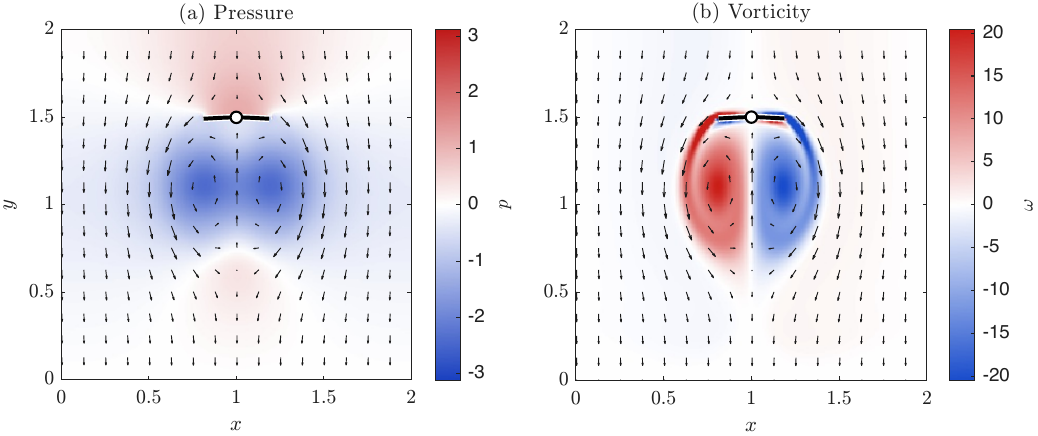}
\caption{Pressure and vorticity at the final time \(t=1.875\) on the \(h=\frac{1}{128}\) grid. Arrows indicate the velocity field in both panels. The black curve shows the fiber configuration, and the circle marks the clamped midpoint.}
\label{fig:asz_channel_snapshot}
\end{figure}

The computations use $N=16$, $32$, $64$, and $128$, with $h=1/N$, on a domain discretized using $2N\times2N$ Eulerian grid cells. We set $\Delta s=h/2$ and $\Delta t=0.05h$. This small Courant number is chosen to avoid instabilities associated with the forward Euler treatment of the convective term. At lower Reynolds numbers, a Courant number closer to unity remains stable with the same bending, stretching, and constraint stiffness parameters. This indicates that the time-step restriction is induced by the convective term rather than by the stiffness of the fiber. The four discretizations employ $13$, $25$, $49$, and $97$ fiber markers, respectively. Each calculation advances from rest to $t=1.875$, corresponding to five convective time units. \Cref{fig:asz_channel_snapshot} shows the flow and fiber configuration at the final time for the finest grid.

\newpage
\setlength{\textfloatsep}{10pt}
\begin{table}[!t]
\centering
\scriptsize
\caption{CAV-preconditioned Newton--Krylov performance for the flexible fiber model under simultaneous Eulerian and Lagrangian refinement. Each calculation uses $\ds=h/2$ and advances to $t=1.875$. Newton counts are averaged over time steps and FGMRES counts are averaged over the individual Newton linear systems. The reported average patch size is computed over the patches on the finest multigrid level. The maximum patch is the largest CAV patch observed over the duration of the simulation.}
\label{tab:asz_refinement}
\begin{adjustbox}{max width=\textwidth}
\begin{tabular}{cccccccc}
\toprule
$N$ & \shortstack{Eulerian\\grid} & $M$ &
\shortstack{Avg. Newton\\its. per step} &
\shortstack{Avg. FGMRES\\its. per solve} &
\shortstack{Max. FGMRES\\its. per solve} &
\shortstack{Avg. patch\\size} &
\shortstack{Max. patch\\size} \\
\midrule
$16$  & $32\times32$   & $13$ & $2.00$ & $8.65$ & $13$ & $11.76$ & $201$ \\
$32$  & $64\times64$   & $25$ & $1.48$ & $9.02$ & $14$ & $7.85$  & $198$ \\
$64$  & $128\times128$ & $49$ & $1.45$ & $9.15$ & $15$ & $6.36$  & $195$ \\
$128$ & $256\times256$ & $97$ & $1.60$ & $9.53$ & $16$ & $5.68$  & $197$ \\
\bottomrule
\end{tabular}
\end{adjustbox}
\end{table}

\Cref{tab:asz_refinement} shows little growth in the linear iteration counts as the number of Eulerian unknowns increases. The average FGMRES count increases from $8.65$ to $9.53$ and no linear solve requires more than $16$ iterations. Newton's method requires at most two corrections per time step at every resolution. The largest CAV patch contains between $195$ and $201$ unknowns throughout the refinement study. Because the IB coupling is confined to a fixed number of grid cells surrounding the fiber, the fraction of pressure cells whose patches are further enlarged decreases under grid refinement, explaining the decrease in the mean patch size observed as the grid is further refined. These results support the bounded-patch argument in Subsection~\ref{subsec:patch_relaxation} and demonstrate essentially grid-independent convergence for the nonlinear flow past a flexible fiber benchmark. The constraint violations also decrease under grid refinement. For \(N=16\), \(32\), \(64\), and \(128\), respectively, the maximum midpoint displacements are \(7.56\times10^{-3}\), \(1.98\times10^{-3}\), \(8.82\times10^{-4}\), and \(4.12\times10^{-4}\). The corresponding maximum midpoint-tangent constraint errors are \(3.08\times10^{-2}\), \(6.22\times10^{-4}\), \(2.34\times10^{-4}\), and \(9.33\times10^{-5}\). The maximum recorded relative errors in the total fiber length are \(0.312\%\), \(0.0535\%\), \(0.0112\%\), and \(0.00296\%\).

We next examine parameter robustness at fixed Eulerian grid resolution of $h = \frac{1}{64}$, corresponding to a $128\times128$ Eulerian grid. Every calculation starts from rest and advances to $tU/L=2$ with $\Delta t/h=0.05$. The baseline parameters are the ones for the grid refinement study above. Let $\gamma_{\mathrm b}$ and $\gamma_{\mathrm s}$ denote multiplicative factors applied to $EI$ and $K_{\mathrm s}$, respectively, and let $\gamma_{\mathrm c}$ denote a common factor applied to $K_{\mathrm p}$ and $K_{\mathrm t}$. Thus, the baseline parameters correspond to taking $\gamma_{\mathrm b}=\gamma_{\mathrm s}=\gamma_{\mathrm c}=1$. We vary the Reynolds number and each stiffness family separately, then combine $\mathrm{Re}=2000$ with a tenfold increase in all four stiffness coefficients to test the robustness of the preconditioner.

\begin{table}[H]
\centering
\scriptsize
\caption{CAV-preconditioned Newton--Krylov performance for flow past a flexible fiber under variation of the Reynolds number and structural stiffnesses at the fixed grid resolution $h = \frac{1}{64}$. Each calculation advances to $tU/L=2$. The first row gives the baseline. The next pairs of rows vary the Reynolds number, all stiffnesses, the bending stiffness, and the stretching stiffness, respectively. The final row combines the largest Reynolds number with the largest simultaneous stiffness scaling. Newton counts are averaged over time steps, and FGMRES counts are averaged over the individual Newton linear systems.}
\label{tab:asz_parameter_robustness}
\begin{adjustbox}{max width=\textwidth}
\begin{tabular}{cccccccc}
\toprule
$\mathrm{Re}$ & $\gamma_{\mathrm b}$ & $\gamma_{\mathrm s}$ & $\gamma_{\mathrm c}$ &
\shortstack{Avg. Newton\\its. per step} &
\shortstack{Avg. FGMRES\\its. per solve} &
\shortstack{Max. FGMRES\\its. per solve} &
\shortstack{Max. patch\\size} \\
\midrule
$200$  & $1$   & $1$   & $1$   & $1.79$ & $9.61$ & $15$ & $195$ \\
\addlinespace[2pt]
$20$   & $1$   & $1$   & $1$   & $1.78$ & $9.84$ & $15$ & $198$ \\
$2000$ & $1$   & $1$   & $1$   & $1.92$ & $9.62$ & $15$ & $195$ \\
\addlinespace[2pt]
$200$  & $0.1$ & $0.1$ & $0.1$ & $1.92$ & $9.51$ & $15$ & $247$ \\
$200$  & $10$  & $10$  & $10$  & $1.69$ & $9.75$ & $15$ & $185$ \\
\addlinespace[2pt]
$200$  & $0.1$ & $1$   & $1$   & $2.00$ & $9.31$ & $15$ & $253$ \\
$200$  & $10$  & $1$   & $1$   & $1.45$ & $9.40$ & $15$ & $198$ \\
\addlinespace[2pt]
$200$  & $1$   & $0.1$ & $1$   & $1.53$ & $9.33$ & $15$ & $198$ \\
$200$  & $1$   & $10$  & $1$   & $1.97$ & $9.37$ & $15$ & $195$ \\
\midrule
$2000$ & $10$ & $10$ & $10$ & $1.88$ & $9.79$ & $15$ & $211$ \\
\bottomrule
\end{tabular}
\end{adjustbox}
\end{table}

\Cref{tab:asz_parameter_robustness} shows that the average FGMRES count remains between $9.31$ and $9.84$ as the Reynolds number ranges from $20$ to $2000$ and the bending, stretching, and clamp stiffnesses are varied by factors of ten in either direction. The maximum number of FGMRES iterations for each study is $15$. Most notably, increasing the Reynolds number from $200$ to $2000$ while increasing all four stiffness coefficients by a factor of ten changes the average FGMRES count only from $9.61$ to $9.79$. This demonstrates that the preconditioner remains effective when weak viscous coupling and strong elastic coupling occur simultaneously.

\section{Conclusions}
\label{sec:conclusions}

We developed and evaluated a coupling-aware Vanka (CAV) geometric multigrid preconditioner for the reduced Eulerian velocity--pressure systems arising from semi-implicit IB discretizations. CAV uses the nonzero graph of the Eulerian elasticity matrix to combine standard pressure-centered Vanka patches, thereby capturing elasticity-induced velocity coupling and retaining the Cartesian-grid structure of the reduced system. Because the spreading, interpolation, and Lagrangian force operators have bounded stencil sizes, the CAV patch sizes remain bounded under grid refinement and each multigrid cycle has linear complexity in the number of Eulerian unknowns.

The two-grid membrane tests provide spectral evidence for the effectiveness of this construction across broad ranges of structural stiffness and fluid viscosity. Over $10^2\le\kappa_{\mathrm s}\le10^6$, every nonzero eigenvalue of the preconditioned two-grid operator has modulus at least $0.809$. In the low-viscosity membrane regime examined here, CAV produces a substantially smaller two-grid spectral radius than $\mathsf B(4,2)$. The simplified additive Schwarz momentum-block analysis in Appendix~\ref{sec:appendix} identifies a complementary mechanism: velocities with zero or small Eulerian elastic energy should admit patchwise decompositions that do not create excess elastic energy in the local components. Together, the spectral results and the Schwarz analysis provide complementary empirical and analytical evidence that CAV captures the stiffness-dominant coupling missed by fixed geometric boxes.

Across the target-point, membrane, and beam tests, CAV-preconditioned FGMRES reduces the relative residual by ten orders of magnitude in $9$--$15$ iterations, with only weak growth under grid refinement and no change to the patch-construction algorithm. The $\mathsf B(4,2)$ iteration counts grow under refinement for the membrane and beam tests but remain competitive for the target-point test, whose Lagrangian operator has no intermarker coupling. Thus, coupling-aware patches are most beneficial for force laws that produce extended coupling along the immersed structure, but the same construction remains effective for a localized force law. The flow past a flexible fiber provides a more demanding test involving stretching, bending, positional, and tangential forces. As the Eulerian grid is refined from $32\times32$ to $256\times256$ cells, the average FGMRES count grows only slightly, from $8.65$ to $9.53$, Newton's method requires at most two corrections per time step, and the maximum CAV patch size remains bounded independent of grid size. Additionally, at a fixed grid resolution, the average FGMRES count remains below $9.8$ even when the Reynolds number and structural stiffness parameters are simultaneously increased.

A rigorous analysis of the pressure-coupled CAV V-cycle remains an important direction for future work. Existing Schwarz, Vanka, and constrained-multigrid theories for Stokes problems provide possible starting points~\cite{schoberl2003saddle,manservisi2006vanka,chen2015constrained}. The candidate weighted additive construction proposed in Equation~\eqref{eq:additive_cav_saddle} also provides a natural route to replacing the sequential multiplicative sweep and improving parallel scalability, provided that robust weights can be determined for the indefinite saddle-point system. Taken together, the present results demonstrate that elasticity-aware combinations of standard Vanka patches provide a robust, grid-scalable geometric multigrid preconditioner with linear per-cycle complexity for time-dependent implicit IB formulations.

\appendix
\section{Additive Schwarz motivation for CAV patch selection}
\label{sec:appendix}

We develop an additive Schwarz model of the elasticity-augmented momentum block to motivate how the CAV patches defined in Subsection~\ref{subsec:cav_patches} capture the elastic near-kernel at high structural stiffness. The model associates a velocity subspace with the velocity unknowns in each CAV patch, derives a stable-decomposition condition for stiffness-robust preconditioning, and motivates a pressure-coupled additive extension. The implemented multigrid method also includes pressure coupling, multiplicative patch ordering, and coarse-grid correction, so the model provides a heuristic interpretation of the full solver.

Assume that the Lagrangian elasticity matrix has the form \(\Eb=\kappa\Ezero\), in which \(\Ezero\) is symmetric negative semidefinite. For uniform Lagrangian spacing \(\ds\) and quadrature weights \(\omega_k=\ds\), the discrete adjoint identity stated in \Cref{sec:eulerian_schur} has the coefficient form \(\Sb=(\ds/h^2)\Jb\tran\). The momentum block can then be written as
\begin{equation*}
\Ab_\kappa=\Azero-\dt\,\kappa\Eeulzero,\qquad
\Azero=\frac{\rho}{\dt}\Ib-\mu\Lb,\qquad
\Eeulzero=-\Cb\tran\Cb,
\end{equation*}
in which
\begin{equation*}
\Cb=\left(\frac{\ds}{h^2}\right)^{1/2}(-\Ezero)^{1/2}\Jb .
\end{equation*}
The matrices \(\Azero\) and \(\Ab_\kappa\) are symmetric positive definite, and \(\Eeulzero\) is symmetric negative semidefinite. Thus, \(\Eeul=\kappa\Eeulzero\) in Equation~\eqref{eq:intro_saddle}. A velocity \(\ub\in\ker\Cb\) has zero elastic energy since \(-\ub\tran\Eeulzero\ub=\|\Cb\ub\|^2=0\).

Recall that the global discrete Eulerian velocity space is \(\Uh\). Let \(\mathscr U_i\subset\Uh\) be the subspace supported on the velocity unknowns in CAV patch \(i\). These spaces satisfy \(\Uh=\sum_i\mathscr U_i\). The Boolean matrix \(\Rb_i\) restricts global velocity coefficient vectors to \(\mathscr U_i\). This patch restriction is distinct from the intergrid restriction \(\Rb_{\ell}^{\ell-1}\) used in the multigrid algorithm of \Cref{sec:multigrid}. The exact additive Schwarz matrix is
\begin{equation}
\label{eq:additive_preconditioner}
\Bb_\kappa
=\sum_i\Rb_i\tran
\left(\Rb_i\Ab_\kappa\Rb_i\tran\right)^{-1}\Rb_i .
\end{equation}
Equation~\eqref{eq:additive_preconditioner} is itself the unweighted additive CAV momentum correction: it uses the velocity portions of the CAV patches defined in Subsection~\ref{subsec:cav_patches}, but computes every local correction from the same residual. Thus, the calculation is not tied to the multiplicative ordering used in the experiments.

Each restricted matrix is positive definite because \(\Ab_\kappa\) is positive definite and \(\Rb_i\) is assumed to be full rank. Setting \(\Mb_\kappa=\Bb_\kappa^{-1}\), the standard additive Schwarz minimum-decomposition identity gives~\cite{xu1992,toselli2005}
\begin{equation}
\label{eq:min_decomp}
\ub\tran\Mb_\kappa\ub
=\min_{\substack{\ub=\sum_i\ub_i\\\ub_i\in\mathscr U_i}}
\sum_i\left(\ub_i\tran\Azero\ub_i
+\dt\,\kappa\|\Cb\ub_i\|^2\right).
\end{equation}
Furthermore, the eigenvalues of \(\Bb_\kappa\Ab_\kappa\) are the stationary values of the generalized Rayleigh quotient
\begin{equation*}
\mathcal Q_\kappa[\ub]
=\frac{\ub\tran\Azero\ub+\dt\,\kappa\|\Cb\ub\|^2}
{\ub\tran\Mb_\kappa\ub}.
\end{equation*}
For a global velocity \(\ub\in\ker\Cb\) and a patchwise decomposition \(\ub=\sum_i\ub_i\), the assembled elastic response satisfies \(\sum_i\Cb\ub_i=\Cb\ub=\zerob\), although the individual local responses \(\Cb\ub_i\) need not vanish. Equation~\eqref{eq:min_decomp} shows that if such a velocity cannot be decomposed as
\begin{equation}
\label{eq:kernel_decomposition}
\ub=\sum_i\ub_i,
\qquad \ub_i\in\mathscr U_i\cap\ker\Cb,
\end{equation}
then every patchwise decomposition will contribute positive elastic energy and \(\mathcal Q_\kappa[\ub]\) can decrease proportionally to \((\dt\,\kappa)^{-1}\). This is the canonical obstruction in subspace correction methods for nearly singular systems~\cite{lee2007nearlysingular,wu2014parallel}.

Conversely, suppose that every \(\ub\in\Uh\) admits a patchwise decomposition satisfying
\begin{align*}
\ub&=\sum_i\ub_i,
&\ub_i&\in\mathscr U_i, \\
\sum_i\ub_i\tran\Azero\ub_i&\le C_{\text{fluid}}\ub\tran\Azero\ub,
&\sum_i\|\Cb\ub_i\|^2
&\le C_{\text{elastic}}\|\Cb\ub\|^2.
\end{align*}
with constants $C_{\text{fluid}}$ and $C_{\text{elastic}}$ independent of $\kappa$. The first inequality controls the total fluid energy of the local components, and the second prevents the decomposition from creating excess elastic energy. Using this decomposition in the minimum on the right-hand side of Equation~\eqref{eq:min_decomp} gives
\begin{equation*}
\ub\tran\Mb_\kappa\ub
\le \max(C_{\text{fluid}},C_{\text{elastic}})
\left(\ub\tran\Azero\ub+\dt\,\kappa\|\Cb\ub\|^2\right),
\end{equation*}
and therefore \(\mathcal Q_\kappa[\ub]\ge1/\max(C_{\text{fluid}},C_{\text{elastic}})\) independently of \(\kappa\). To bound the other end of the spectrum, patch spaces \(\mathscr U_i\) and \(\mathscr U_j\) interact if
\begin{equation*}
\Rb_i\Azero\Rb_j\tran\ne\zeromat
\qquad\text{or}\qquad
\Rb_i\Eeulzero\Rb_j\tran\ne\zeromat.
\end{equation*}
The quantity \(N_{\text{int}}\) is the maximum number of patch spaces interacting with any one \(\mathscr U_i\), including \(\mathscr U_i\) itself. Thus, interaction includes not only overlapping patches, but also disjoint patches whose degrees of freedom are coupled by \(\Azero\) or \(\Eeulzero\). For any patchwise decomposition \(\ub=\sum_i\ub_i\), the Cauchy--Schwarz inequality in the \(\Ab_\kappa\) inner product gives
\begin{align*}
\ub\tran\Ab_\kappa\ub
&=\sum_i\sum_{j\in\mathcal I_i}
  \ub_i\tran\Ab_\kappa\ub_j \\
&\le \sum_i\sum_{j\in\mathcal I_i}
  \|\ub_i\|_{\Ab_\kappa}\|\ub_j\|_{\Ab_\kappa} \\
&\le N_{\text{int}}\sum_i\|\ub_i\|_{\Ab_\kappa}^2,
\end{align*}
in which \(\mathcal I_i\) is the set of patches interacting with \(\mathscr U_i\) and the last inequality follows from \(2ab\le a^2+b^2\). Here, \(\|\vb\|_{\Ab_\kappa}^2=\vb\tran\Ab_\kappa\vb\). Applying this estimate to a decomposition that attains the minimum in Equation~\eqref{eq:min_decomp} yields \(\mathcal Q_\kappa[\ub]\le N_{\text{int}}\). Hence, the spectrum of the resulting preconditioned matrix lies in \(\left[1/\max(C_{\text{fluid}},C_{\text{elastic}}),N_{\text{int}}\right]\).
In typical IB discretizations, the regularized delta kernel spans a fixed
number of grid cells, and the constitutive force is approximated by a
fixed-width Lagrangian stencil with \(\ds=\mathcal O(h)\). Each CAV patch
thus interacts with only $\mathcal O(1)$ other patches, so $N_{\text{int}}$
is bounded independently of the grid resolution and the scalar model
parameters.

CAV addresses the lower end of the spectrum by placing velocities coupled through \(\Eeulzero\) in the same patch, giving each local space the degrees of freedom needed to represent components with zero or small elastic energy. This construction is intended to enforce the kernel-recovery condition in Equation~\eqref{eq:kernel_decomposition} and keep the stable-decomposition constants $C_{\text{fluid}}$ and $C_{\text{elastic}}$ bounded under mesh refinement and parameter variation. Under these conditions, the additive momentum-block spectrum remains in a fixed interval, and the corresponding Krylov convergence is robust. The implemented method also includes pressure coupling, multiplicative patch ordering, and coarse-grid correction, so the additive result serves only as a qualitative guide. The empirical spectra in \Cref{subsec:membrane} assess this prediction for the implemented pressure-coupled cycle.

The candidate pressure-coupled additive relaxation uses \(\widetilde{\Rb}_i\), the extension of the velocity restriction \(\Rb_i\) to the full velocity--pressure CAV patch, and diagonal local weights \(\Thetab_i\) satisfying the partition-of-unity condition
\begin{equation*}
\sum_i\widetilde{\Rb}_i\tran\Thetab_i\widetilde{\Rb}_i=\Ib.
\end{equation*}
The candidate operator is
\begin{equation}
\label{eq:additive_cav_saddle}
\Bb_{\text{CAV}}^{\text{add}}
=
\alpha\sum_i\widetilde{\Rb}_i\tran\Thetab_i
\left(\widetilde{\Rb}_i\KE
\widetilde{\Rb}_i\tran\right)^{-1}
\widetilde{\Rb}_i .
\end{equation}
Here, \(\alpha\) is a relaxation factor. All multiplicative tests reported in this paper use \(\alpha=1\). A simple choice weights each unknown by the reciprocal of its patch count, which yields a partition of unity. Restricted additive injection is another option. The heuristic momentum-block model motivates the patch family and suggests which near-kernel components it should capture, but damping, weighting, and the indefinite pressure block require additional analysis. Equation~\eqref{eq:additive_cav_saddle} therefore provides a concrete candidate for pressure-coupled additive CAV relaxation; its robustness for the indefinite system remains to be established.

\section*{Code and data availability}
The code needed to reproduce the experiments is available at
\url{https://github.com/IBAMR/coupling-aware-vanka-reproducibility}.
 
\section*{Acknowledgments}
\TheFunding

\renewcommand\refname{REFERENCES}

\end{document}